
%
%

%
%

%
%

\font\scriptsize=cmr7

\input amstex

\documentstyle{amsppt}
  \magnification=1000
  \hsize=7.0truein
  \vsize=9.0truein
  \hoffset -0.1truein
  \parindent=2em

\define\Ac{{\Cal A}}                           

\define\Ao{\overset {\scriptsize o}\to A}      


\define\Bc{{\Cal B}}                           

\define\Bo{\overset {\scriptsize o}\to B}      

\define\Cc{{\Cal C}}                           

\define\Ccw{\widetilde{\Cal C}}                

\define\Cpx{\bold C}                           

\define\dif{\text{\it d}}                      

\define\dcup#1{\underset#1\to\cup}             

\define\freeF{\bold F}                         

\define\Intg{\bold Z}                          

\define\lspan{\text{\rm span}@,@,@,}           

\define\MvN{{\Cal M}}                          

\define\Nats{\bold N}                          

\define\NvN{{\Cal N}}                          

\define\NvNo{\overset{\scriptsize o}\to \NvN}  

\define\NvNw{\widetilde{\NvN}}                 

\define\pol{\text{\rm pol}}                    

\define\pw{\tilde{p}}                          

\define\QED{$\hfill$\qed\enddemo}              

\define\Qrat{\bold Q}                          

\define\Reals{\bold R}                         

\define\restrict{\lower .3ex                   
     \hbox{\text{$|$}}}

\define\Rw{\widetilde{R}}                      

\topmatter

  \title Interpolated free group factors \endtitle

  \author Ken Dykema \endauthor

  \affil University of California, \\
         Berkeley, California, USA 94720, \\
         (e--mail dykema\@math.berkeley.edu) \endaffil

  \thanks Studies and research
          supported by the Fannie and John Hertz Foundation. \endthanks
  \thanks This work will form part of the author's Ph.D\. thesis
          at the University of California, Berkeley. \endthanks

  \date March 1992 \enddate

  \abstract  The interpolated free
             group factors $L(\freeF_r)$  for $1<r\le\infty$, (also defined
             by F\. R\u{a}dulescu) are given another (but equivalent)
             definition as well as proofs
             of their properties with respect to compression by projections
             and free products.  In order to prove the addition formula for
             free products,
             algebraic techniques are developed
             which allow us to show $R*R\cong L(\freeF_2)$ where $R$ is the
             hyperfinite II$_1$--factor.
  \endabstract

\endtopmatter

\document
  \TagsOnRight
  \baselineskip=14pt

\noindent{\bf Introduction.}

  The free group factors $L(\freeF_n)$ for $n=2,3,\ldots,\infty$
    (introduced in~\cite{4}) have recently been extensively
    studied~\cite{11,2,5,6,7}
    using Voiculescu's theory of freeness in noncommutative probability
    spaces (see~\cite{8,9,10,11,12,13}, especially the latter for
    an overview).
  One hopes to eventually be able to solve the old
    isomorphism question,
    first raised by R.V\. Kadison in the 1960's, of whether
    $L(\freeF_n)\cong L(\freeF_m)$ for $n\neq m$.
  In~\cite{7}, F\. R\u{a}dulescu introduced II$_1$--factors $L(\freeF_r)$
    for $1<r\le\infty$, equaling the free group factor $L(\freeF_n)$
    when $r=n\in\Nats\backslash\{0,1\}$ and satisfying
    $$ L(\freeF_r)*L(\freeF_{r'})=L(\freeF_{r+r'}),\;(1<r,r'\le\infty)
      \tag1 $$
    and
    $$ L(\freeF_r)_\gamma=L(\freeF(1+\frac{r-1}{\gamma^2})),\;
      (1<r\le\infty,\,0<\gamma<\infty). \tag2 $$
  Where for a II$_1$--factor $\MvN$, $\MvN_\gamma$ means the
    algebra~\cite{4} defined as follows:
    for $0<\gamma\le1$, $\MvN_\gamma=p\MvN p$, where $p\in\MvN$
      is a self--adjoint projection of trace $\gamma$;
     for $\gamma=n=2,3,\ldots$ one has $\MvN_\gamma=\MvN\otimes M_n(\Cpx)$;
     for $0<\gamma_1,\gamma_2<\infty$ one has
    $$\MvN_{\gamma_1\gamma_2}= (\MvN_{\gamma_1})_{\gamma_2}. $$

  We had independently found the interpolated free group factors $L(\freeF_r)$
    $(1<r\le\infty)$ and the formulas~(1) and~(2),
    defining them
    differently and using different techniques.
  In this paper we give our definition and proofs.
  This picture of $L(\freeF_r)$ is sometimes more convenient,
    {\it e.g\.}~\S4 of~\cite{3}.
  It is a natural extension of the result~\cite{2} that
    $$ \L(\Intg)*R\cong\L(\freeF_2), \tag3 $$
    where $R$ is the hyperfinite
    II$_1$--factor.
  We also introduce some elementary algebraic techniques for freeness which
    have further application in~\cite{3}.
  One consequence of them that we prove here is that $R*R\cong\L(\freeF_2)$.

  This paper has four sections.  In~\S1 we state a random matrix result
    (from~\cite{2}, \cite{7}) and some consequences;
   in~\S2 we define the interpolated free group factors and prove the
    formula~(2);
   in~\S3 we develop the algebraic techniques;
   in~\S4 we prove
    the addition formula~(1) and also make an observation
    from~(1) and~(2) (also observed in~\cite{7}) that,
    as regards the isomorphism question, we must have one of
   two extremes.
  Our original proof of the addition formula~(1) was a fairly messy
    application
    of the algebraic techniques developed in~\S3.
  The proof of Theorem~4.1
    that appears here, while still using the algebraic techniques
    in an essential way, benefits significantly from ideas found in
    the proof of
    F\.~R\u adulescu~\cite{7}.

\noindent{\bf \S1\. The matrix model.}

  Voiculescu, as well as developing the whole notion of freeness in
    noncommutative probability spaces, had the fundamental idea of
    using Gaussian random matrices to model freeness, which he developed
    in~\cite{12}.
  In~\cite{2}, we extended this matrix model to the non--Gaussian case
    and also to be able to handle
    semicircular families together with a free finite dimensional algebra.
  As R\u{a}dulescu observed in~\cite{7}, the
    matrix model necessary to be able to handle the free
    finite dimensional algebra can be easily proved in the
    Gaussian case directly using Voiculescu's methods
    ({\it cf\.} the appendix of~\cite{2}).
  In any case, we shall use this matrix model in this paper, and quote it
    here, as well as some results of it.
  Our notation for random matrices will be as in~\cite{2}.
  A trivial reformulation of Theorem~2.1 of~\cite{2} gives
  \proclaim{Theorem 1.1}
    Let $Y(s,n)\in M_n(L)$ for $s\in S$ be self--adjoint
      independently distributed
      $n\times n$ random matrices as in Theorem~2.1 of~\cite{2}.
    For $c=\left(\matrix c_{11} & \cdots & c_{1N} \\
                         \vdots & \ddots & \vdots \\
                         c_{N1} & \cdots & c_{NN}
      \endmatrix \right) \in M_N(\Cpx)$
      and for $n$ a multiple of $N$ let
      $c(n)=\left(\matrix
         c_{11}I_{\frac nN} & \cdots & c_{1N}I_{\frac nN} \\
         \vdots            & \ddots & \vdots \\
         c_{N1}I_{\frac nN} & \cdots & c_{NN}I_{\frac nN}
      \endmatrix\right)$
      be a constant matrix in $M_n(L)$.
    Then
      $$ \bigl\{(\{Y(s,n)\})_{s\in S},\{c(n)\mid c\in M_N(\Cpx)\}\bigr\} $$
      is an asymptotically free family as $n\rightarrow\infty$, and each
      $Y(s,n)$ has for limit distribution a semicircle law.
  \endproclaim

  An immediate result of the above is (3.2 of~\cite{2})
  \proclaim{Theorem 1.2}
    In a noncommutative probability space $(\MvN,\phi)$ with $\phi$
      a trace, let $\nu_1=\{X^s\mid s\in S\}$ be a semicircular family
      and let $\nu_2=\{e_{ij}\mid1\le i,j\le n\}$ be a system of matrix
      units such that $\{\nu_1,\nu_2\}$ is free.
    Then in $(e_{11}\MvN e_{11},n\phi\restrict_{e_{11}\MvN e_{11}})$,
      $\omega_1=\{e_{1i}X^se_{i1}\mid1\le i\le n,\,s\in S\}$ is a
      semicircular family and
      $\omega_2=\{e_{1i}X^se_{j1}\mid1\le i<j\le n,\,s\in S\}$ is a
      circular family such that $\{\omega_1,\omega_2\}$ is free.
  \endproclaim

  The following is analogous to Theorem~2.4 of~\cite{11}.
  \proclaim{Theorem 1.3}
    In a noncommutative probability space $(\MvN,\phi)$ with $\phi$
      a trace, let $\nu=\{X^s\mid s\in S\}$ be a semicircular family
      and let $R$ be a copy of the hyperfinite II$_1$--factor such that
      $\{\nu,R\}$ is free.
    Let $p\in R$ be a nonzero self--adjoint projection.
    Then in $(p\MvN p,\phi(p)^{-1}\phi\restrict_{p\MvN p})$,
      $\omega=\{pX^sp\mid s\in S\}$ is a semicircular family and
      $\{pRp,\omega\}$ is free.
    (Note from~\cite{4} that $pRp$ is also a copy of the hyperfinite
      II$_1$--factor.)
  \endproclaim
  \demo{Proof}
    Suppose first that $\phi(p)=m/2^k$, a dyadic rational number.
    Since for $U\in R$ a unitary, $\{R,U\nu U^*\}$ is free, we may
      let $p$ be any projection in $R$ of the given trace.
    Writing $R=M_{2^k}\otimes M_2\otimes M_2\otimes\cdots$,
      we use Theorem~1.1 in order to
      model $\nu$ as the limit of self-adjoint independently distributed
      random matrices of size $n=2^k,2^{k+1},2^{k+2},\ldots$,
      and model a dense subalgebra of $R$ (equal to the tensor product
      of matrix algebras)
      by constant random matrices.
    Choosing $p$ to correspond to a diagonal element of $M_{2^k}$,
      we may apply Theorem~1.1 again to see that $\omega$ is a
      semicircular family, $pRp\cong M_m\otimes M_2\otimes M_2\otimes\cdots$,
      and $\{pRp,\omega\}$ is free.

    Now for general $p$, let $(p_l)_{l=1}^\infty$ be a decreasing
      sequence of projections in $R$ which converge to $p$ and such
      that each $\phi(p_l)$ is a dyadic rational number.
    Then
      $$ \bigl\{p_lRp_l=\{p_lyp_l\mid y\in R\},
        \{p_lX^sp_l\mid s\in S\}\bigr\} $$
      has limit distribution equal to $\{pRp,\omega\}$ as $l\rightarrow\infty$.
    For each $l$ we have freeness and semicircularity,
      hence also in the limit.
  \QED

  In addition, modeling $R$ and a semicircular family
    as in the above proof, we can easily prove
  \proclaim{Theorem 1.4}
    In a noncommutative probability space $(\MvN,\phi)$ with $\phi$
      a trace, let $\nu=\{X^s\mid s\in S\}$ be a semicircular family,
      and $R$ a hyperfinite II$_1$--factor containing a system of matrix
      units $\{e_{ij}\mid1\le i,j\le n\}$,
      such that $\{\nu,R\}$ is free.
    Then in $(e_{11}\MvN e_{11},n\phi\restrict_{e_{11}\MvN e_{11}})$,
      $\omega_1=\{e_{1i}X^se_{i1}\mid1\le i\le n,\,s\in S\}$ is a
      semicircular family and
      $\omega_2=\{e_{1i}X^se_{j1}\mid1\le i<j\le n,\,s\in S\}$ is a
      circular family such that $\{\omega_1,\omega_2,e_{11}Re_{11}\}$
      is $*$--free.
  \endproclaim

\noindent{\bf \S2\. Definition and compressions of $L(\freeF_r)$.}

  \proclaim{Definition 2.1} \rm
    In a W$^*$--probability space $(\MvN,\tau)$, where $\tau$ is a
      faithful trace, let $R$ be a copy of the hyperfinite II$_1$--factor
      and $\omega=\{X^t\mid t\in T\}$ be a semicircular family such that
      $R$ and $\omega$ are free.
    Then $L(\freeF_r)$ for $1<r\le\infty$ will denote
      any factor isomorphic to $(R\cup\{p_t X^t p_t\mid t\in T\})''$,
      where $p_t\in R$ are self--adjoint projections
      and $r=1+\sum_{t\in T}\tau(p_t)^2.$
  \endproclaim

  \proclaim{Proposition 2.2}
    $L(\freeF_r)$ is well--defined, {\rm i.e\.} if
    $\Ac=(R\cup\{p_t X^t p_t\mid t\in T\})''$ and
    $\Bc=(R\cup\{q_t X^t q_t\mid t\in T\})''$, where
    $1+\sum\tau(p_t)^2=r=1+\sum\tau(q_t)^2$, then $\Ac\cong\Bc$.
  \endproclaim
  \demo{Proof}
    We show that $\Ac$ (and thus also $\Bc$) is isomorphic to an algebra
      of a certain ``standard form.''
    Let $(f_k)_{k=1}^\infty$ be an orthogonal family of projections in
      $R$ such that $\tau(f_k)=2^{-k}$, and let $f_0=1$.
    We show that $\Ac$ is isomorphic to
      $\Cc=(R\cup\{f_{k_s}X^sf_{k_s}\mid s\in S\})''$, where $S\subseteq T$,
      each $k_s\in\Nats=\{0,1,2,\ldots\}$, $1+\sum_{s\in S}\tau(f_{k_s})^2=r$
      and for each $k\ge0$, letting $S(k)=\{s\in S\mid k_s>k\}$
      we have that $\sum_{s\in S(k)}2^{-k_s}<2^{-k}$.
    (Note that these conditions imply a unique choice of $|\{s\mid k_s=k\}|$
      for all $k$.)
    This will prove the proposition.

    Proving $\Ac\cong\Cc$ is an exercise in cutting and pasting.
    Note
      that if $U_t$ are unitaries in $R$, $(t\in T)$, then
      $\{R,(\{U_tX^tU_t^*\})_{t\in T}\}$ is free in $(\MvN,\tau)$.
    Moreover, each projection $p\in R$ is conjugate by a unitary in $R$ to a
      projection that is a (possibly infinite) sum of projections
      in $\{f_k\mid k\ge1\}$.
    Hence letting $T'=\{t\in T\mid p_t\neq0\}$,
      we may assume without loss of
      generality that each $p_t$ for $t\in T'$
      is equal to such a sum, and we write
      $p_t=\sum_{k\in K_t}f_k$, for $K_t\subseteq\Nats\backslash\{0\}$
      whenever $t\in T'$ and $p_t\neq1$, and we set $K_t=\{0\}$
      if $p_t=1$.
    Then
      $$ \Ac=\bigl(R\cup
        \{f_kX^tf_{k'}\mid k,k'\in K_t,\,k'\le k,\,t\in T'\}\bigr)''. $$
    Now we may appeal to the matrix model~(\S1)
      to see that (enlarging $T$ if necessary),
      $$ \Ac=\bigl(R\cup\{f_kX^{\alpha(k,k',t)}f_{k'}
        \mid k,k'\in K_t,\,k'\le k,\,t\in T'\}\bigr)'', $$
      where $\alpha$ is a 1--1 map from
      $\{k,k'\in K_t,\,k'\le k,\,t\in T'\}$ onto a subset $T''$ of $T$.
    (The truth of the above assertion is most easily demonstrated when
      $T'$ and each $K_t$ are finite;
     the general case then follows by taking inductive limits.)

    Consider for a moment $f_kX^tf_{k'}$ for $k'<k$, $t\in T$.
    Note that $f_{k'}$ is the sum of $2^{k-k'}$ orthogonal projections,
      each of which is equivalent in $R$ to $f_k$.
    Using the matrix model shows that
      $$ (R\cup\{f_kX^tf_{k'}\})''
        \cong\bigl(R\cup\{f_kX^{t_j}f_k\mid1\le j\le2^{k-k'}\}
        \cup\{f_kX^{t'_j}f_k\mid1\le j\le2^{k-k'}\}\bigr)'', \tag4 $$
      where $t_1,\ldots,t_{2^{k-k'}},t'_1,\ldots,t'_{2^{k-k'}}$
      are distinct elements of $T$, and the isomorphism in~(4)
      maps $R$ identically onto itself.
    Using inductive limits, one obtains
      $$ \Ac\cong\Ccw=
        \bigl(R\cup\{f_{k_s}X^sf_{k_s}\mid s\in S'\}\bigr)'', \tag5 $$
      for $S'$ some subset of $T$, $k_s\in\Nats$ for each $s\in S'$.
    Moreover, checking the arithmetic of the above moves shows that
      $1+\sum_{s\in S'}\tau(f_{k_s})^2=r$.

    Now for the pasting.
    Note that by the matrix model,
      $$ (R\cup\{f_kX^{t_i}f_k\mid1\le i\le4\})''
        \cong (R\cup\{f_{k-1}X^tf_{k-1}\})'' \tag6 $$
      by an isomorphism mapping $R$ identically to itself,
      whenever $k\ge1$,
      $t_1,\ldots,t_4$ are distinct elements of $T$ and $t\in T$.
    Suppose $r<\infty$.
    Now using~(6) and induction, we can show that
      $\Ccw\cong(R\cup\{f_{k_s}X^sf_{k_s}\mid s\in S''\})''$ for some
      $S''\subseteq T$, $k_s\in\Nats$ for each $s\in S''$, and where
      $S''_k=\{s\in S''\mid k_s=k\}$ has cardinality~$0$, $1$, $2$ or~$3$
      whenever $k\in\Nats\backslash\{0\}$.
    (This is a complicated induction argument, which terminates because
      $r<\infty$.)
    We need now be concerned only with the case $|S''_k|=3\;\forall k\ge k'$,
      some fixed $k'>0$.
    In this case, to show that $\Ccw$ is isomorphic to an algebra
      $\Cc$ that is in
      standard form, it will suffice to show
      $$ \bigl(R\cup\{f_{k_s}X^sf_{k_s}\mid s\in\dcup{k\ge k'}S''_k\}\bigr)''
        \cong (R\cup\{f_{k'-1}X^tf_{k'-1}\})'', \tag7 $$
      where $t\in T$ and where the isomorphism maps $R$ identically to
      itself.
    With the matrix model, we easily see that for $N\ge k'$,
      $$ \bigl(R\cup
        \{f_{k_s}X^sf_{k_s}\mid s\in\dcup{k'\le k\le N}S''_k\}\bigr)''
        \cong (R\cup\{f_{k'-1}X^tf_{k'-1}-f_NX^tf_N\})'', $$
      and using inductive limits gives~(7).

    If $r=\infty$, then considering $S'$ from~(5) and letting
      $S'_k=\{s\in S'\mid k_s=k\}$, we have that
      $\sum_{k=0}^\infty|S'_k|4^{-k}=\infty$.
    Now by repeated application of~(6), we can transform the situation
      (by isomorphisms mapping $R$ identically to itself)
      so that first some $|S'_k|=\infty$, then all
      $|S'_k|=\infty$, then $|S'_0|=\infty$ and $|S'_k|=0$ for all $k\ge1$.
    Thus $\Ccw=L(\freeF_\infty)$ by~(3).
  \QED

  \proclaim{Remark 2.3} \rm
    Formula~(2), together with the fact that $L(\freeF_r)$ for $r\in\Nats$
      is the free group factor on $r$ generators, shows that
      Definition~2.1 is equivalent to R\u{a}dulescu's
      definitions~4.1 and~5.3 of~\cite{7}.
    However, for $r\ge2$, ({\it i.e\.} R\u{a}dulescu's~4.1),
      this equivalence can be seen directly using the ``standard form''
      of $L(\freeF_r)$ as defined in Proposition~1.3, and by noting
      that the isomorphism
      $$ R*L(\Intg)\overset\sim\to\rightarrow L(\Intg)*L(\Intg) \tag8 $$
      in~\cite{2} sends the set of projections $\{f_k\mid k\ge1\}\subset R$
      into one of the copies of $L(\Intg)$ on the right hand side of~(8).
  \endproclaim

  The formula in the following theorem
    for the compression of an interpolated free group factor $L(\freeF_r)$
    by a projection of trace $\gamma$
    was first proved by Voiculescu~\cite{11} for the cases
    $r=2,3,\ldots$, $\gamma=\frac12,\frac13,\frac14,\ldots$ and
    $r=\infty$, $\gamma\in\Qrat_+$.
  It was then extended by F\. R\u{a}dulescu in~\cite{5} for
    $r=\infty$ and $\gamma\in\Reals_+$,
    and in~\cite{6} for $r=2,3,\ldots$ and
    $\gamma=\frac1{\sqrt{2}},\frac1{\sqrt{3}},\ldots$.
  Of course, R\u{a}dulescu also proved this theorem in the generality
    stated here in~\cite{7}.

  \proclaim{Theorem 2.4}
    $$ L(\freeF_r)_\gamma=L(\freeF(1+\frac{r-1}{\gamma^2})) \tag9 $$
    for $1<r\le\infty$ and $0<\gamma<\infty$.
  \endproclaim

  \demo{Proof}
    It suffices to show the case $0<\gamma<1$.
    Let $L(\freeF_r)=\Ac=(R\cup\{p_tX^tp_t\mid t\in T\})''$ be as in
      Definition~2.1, so $1+\sum_{t\in T}\tau(p_t)^2=r$.
    Let $p\in R$ be a projection having trace $\gamma$.
    Without loss of generality, we may assume that each $p_t\le p$.
    Then $p\Ac p=(pRp\cup\{p_tX^tp_t\mid t\in T\})''$, which by Theorem~1.3
      is an interpolated free group factor.  Counting gives the formula~(9).
  \QED

\noindent{\bf \S3\. Algebraic techniques.}

  A crucial ingredient of our proof of the addition formula for free
    products~(1) will be showing that $R*R\cong R*L(\Intg)$, with the
    isomorphism being the identity map on the first copy of $R$.
  In order to show this, we will introduce some elementary techniques
    (Definition~3.4, proof of Theorem~3.5) that are algebraic in nature.
  These techniques have extensive further applications to free products,
    as will be seen in~\cite{3}.

  \proclaim{Remark 3.1} \rm
    In this section, all von Neumann algebras will be finite and have
      fixed normalized faithful traces associated to them,
      and all isomorphisms and inclusions of von Neumann algebras
      will be assumed to be trace preserving.
    Von Neumann algebras that we obtain from others by certain
      operations will have associated traces
      given by the following conventions:
      \roster
      \item"(1)" group von Neumann algebras $L(G)$ for $G$
        a discrete group will have their canonical traces (equal to the
        vector--state for the vector $\delta_e\in l^2(G)$);
      \item"(2)" factors, such as matrix algebras $M_n=M_n(\Cpx)$ or
        the hyperfinite II$_1$--factor $R$, will have (of course)
        their unique normalized traces;
      \item"(3)" a tensor product $A\otimes B$ of algebras will have
        the tensor product trace $\tau_A\otimes\tau_B$
        of the given traces on $A$ and $B$;
      \item"(4)" a free product $A*B$ of algebras will have
        the free product trace $\tau_A*\tau_B$
        of the given traces on $A$ and $B$;
      \item"(5)" if $\MvN$ is a von Neumann algebra with faithful trace
        $\tau$, and $p$ is a projection in $\MvN$, then $p\MvN p$ will
        have trace $\tau(p)^{-1}\tau\restrict_{p\MvN p}$.
      \endroster
    Also, if $A$ is an algebra with specified trace, $\Ao$ will denote
      the ensemble of elements of $A$ whose trace is zero.
  \endproclaim

  First we examine $L(\Intg_2)*L(\Intg_2)$, (where $\Intg_2$ is the two
    element group).
  The fact that $\MvN=L(\Intg_2*\Intg_2)\cong L(\Intg)\otimes M_2$ is
    well known, but we will need the following picture of $\MvN$.
  \proclaim{Proposition 3.2}
    Consider $\MvN=L(\Intg_2)*L(\Intg_2)$ with trace $\tau$, and let
      $p$ and $q$ be projections of trace $\frac12$ generating the first
      and respectively the second copy of $L(\Intg_2)$.
    Then
      $$ \MvN\cong L^\infty([0,\tfrac\pi2],\nu)\otimes M_2, \tag10 $$
      where $\nu$ is a probability measure on $[0,\frac\pi2]$
      without atoms and $\tau$ is given by integration with respect to
      $\nu$ tensored with the normalized trace on $M_2=M_2(\Cpx)$.
    Moreover, in the setup of~(10), we have that
      $$ p=\pmatrix 1 & 0 \\ 0 & 0 \endpmatrix \text{ and }
        q=\pmatrix\cos^2\theta & \cos\theta\sin\theta \\
                  \cos\theta\sin\theta & \sin^2\theta \endpmatrix, \tag11 $$
      where $\theta\in[0,\tfrac\pi2]$.
  \endproclaim
  \demo{Proof}
    It is well known that the universal unital C$^*$--algebra
      generated by two projections $p$ and $q$ is
      $A=\{f:[0,\tfrac\pi2]\rightarrow M_2(\Cpx)\mid f(0)\text{ and }
        f(\frac\pi2)\text{ diagonal}\}$,
      with $p$ and $q$ as in~(11).
    $\MvN$ thus has a dense subalgebra equal to a quotient of $A$, and $\tau$
      gives a trace on $A$.
    One can easily see that a trace on $A$ must be of the following form.
    Let $f(t)_1$ and $f(t)_2$ be the diagonal values of $f(t)$
      for $t=0$ or $\tfrac\pi2$.
    Then
      $$ \tau(f)=a_1f(0)_1+a_2f(0)_2+\int_0^{\tfrac\pi2}\tau_2(f(t))\dif\nu(t)
        + b_1f(\tfrac\pi2)_1+b_2f(\tfrac\pi2)_2, $$
      where $\tau_2$ is the normalized trace on $M_2(\Cpx)$,
      $\nu$ is a positive measure on $[0,\tfrac\pi2]$,
      $a_1,a_2,b_1,b_2\ge0$ and $|\nu|+a_1+a_2+b_1+b_2=1$.
    By Example~2.8 of~\cite{9}, the distribution of $pqp$ in $p\MvN p$
      has no atoms, which implies that $|\nu|=1$ and $\nu$ has no atoms.
  \QED
  \proclaim{Remark 3.3} \rm
    In the right hand side of~(10), let
      $$ x=\pmatrix 0 & 0 \\ 1 & 0 \endpmatrix = \pol((1-p)qp), $$
      where ``$\pol$'' means ``polar part of.''
    Then $x$ is a partial isometry from $p$ to $1-p$ and $\MvN$ is
      generated by $pqp$ together with $x$.
    Let $y=\pol((1-q)pq)$.
    Then $y$ is a partial isometry from $q$ to $1-q$.
    Let
      $$ w=\pmatrix \cos\theta & -\sin\theta \\
                    \sin\theta & \cos\theta \endpmatrix. $$
    Then $w$ is  unitary and $wpw^*=q$, $wxw^*=y$.
  \endproclaim

  \proclaim{Definition 3.4} \rm
    Let $(S_\iota)_{\iota\in I}$ be subsets of a unital algebra $A\ni1$.
    A nontrivial {\it traveling product} in $(S_\iota)_{\iota\in I}$
      is a product $a_1a_2\cdots a_n$ such that $a_j\in S_{\iota_j}$
      $(1\le j\le n)$ and $\iota_1\neq\iota_2\neq\iota_3\neq\cdots\neq\iota_n$.
    The {\it trivial traveling product} is the identity
      element~$1$.
    $\Lambda((S_\iota)_{\iota\in I})$ denotes the set of all traveling
      products in $(S_\iota)_{\iota\in I}$, including the trivial one.
    If $|I|=2$, we will often call traveling products
      {\it alternating products}.
  \endproclaim

  \proclaim{Theorem 3.5}
    Let $A$ and $B$ be finite von Neumann algebras (with specified faithful
      traces---see Remark~3.1).
    Then
      $$ \matrix \format\l&\r&\l\;&\l\\
        \text{(i)}    & (A\otimes L(\Intg_2))&*(B\otimes L(\Intg_2))&
          \cong\; (A*A*B*B*L(\Intg))\otimes M_2, \\ \vspace{1\jot}
        \text{(ii)}   & (A\otimes M_2       )&*(B\otimes L(\Intg_2))&
          \cong\; (A*B*B*L(\freeF_2))\otimes M_2, \\ \vspace{1\jot}
        \text{(iii) } & (A\otimes M_2       )&*(B\otimes M_2       )&
          \cong\; (A*B*L(\freeF_3))\otimes M_2.
      \endmatrix $$
  \endproclaim
  \demo{Proof}
    Let $\MvN$ be the von Neumann algebra on the left hand side of~(i)
      with trace $\tau$.
    It will be notationally convenient to identify $A$ with
      $A\otimes1\subseteq\MvN$ and $B$ with $B\otimes1\subseteq\MvN$.
    Let $p$ and $q$ be projections of trace $\frac12$ contained in the copy
      of $1\otimes L(\Intg_2)$ that commute with $A$ and respectively $B$.
    Let $\NvN_0=\{p,q\}''\cong L(\Intg_2)*L(\Intg_2)$, and let $x,y,w\in\NvN_0$
      be as in Remark~3.3.
    Then
      $$ p\MvN p=
        \bigl(\{pqp\}\cup pA\cup x^*Ax\cup w^*qBw\cup w^*y^*Byw\bigr)''. $$
    We claim moreover that $\{\{pqp\},pA,x^*Ax,w^*qBw,w^*y^*Byw\}$ is
      a free family in $p\MvN p$, which then clearly implies~(i).

    Let us first show that $\{\{pqp\},pA,x^*Ax\}$ is free in $p\MvN p$.
    Let $g_k=(pqp)^k-2\tau((pqp)^k)p$, $(k\ge1)$.
    To show freeness means to show that a nontrivial traveling product
      in $\{g_k\mid k\ge1\}$, $p\Ao$ and $x^*\Ao x$ has trace zero.
    Regrouping gives a traveling product in
      $\Omega_0=\{x,x^*\}\cup\{g_k,xg_k,g_kx^*,xg_kx^*\mid k\ge1\}$
      and $\Ao$.
    Let $a=p-\frac12$, $b=q-\frac12$.
    Then $\NvN_0=\{a,b\}''$, and $\lspan\Lambda(\{a\},\{b\})$ is a dense
      $*$--subalgebra of $\NvN_0$.
    Note that $\Omega_o\subset\NvN_0$, so that by the Kaplansky Density
Theorem,
      any $z\in\Omega_0$ is the s.o.--limit of a bounded sequence in
      $\lspan\Lambda(\{a\},\{b\})$.
    Note also that since $a$ and $b$ are free and each has trace zero, the
      trace of an element of $\lspan\Lambda(\{a\},\{b\})$ is equal to the
      coefficient of $1$.
    Since $\tau(z)=0$, we may chose that approximating sequence
      in $\lspan\Lambda(\{a\},\{b\})$ so that each coefficient of $1$
      equals zero.
    Moreover, since also $\tau(pz)=0$, we may also insist that each
      coefficient of $a$ be zero, {\it i.e\.} we have a bounded approximating
      sequence for $z$ of elements of
      $\lspan(\Lambda(\{a\},\{b\})\backslash\{1,a\})$.
    We must now only show that a nontrivial alternating product in
      $\Lambda(\{a\},\{b\})\backslash\{1,a\})$ and $\Ao$ has trace zero.
    Regrouping gives a nontrivial alternating product in
      $\{a\}\cup\Ao\cup a\Ao$ and $\{b\}$, which by freeness has trace zero.

    Let $\NvN_1=(A\cup\NvN_0)''$, and let us show that
      $\{qw\NvN_1w^*,qB,y^*By\}$ is free in $q\MvN q$, which will complete
      the proof of~(i).
    We show that a nontrivial traveling product in $w\NvN_1w^*$, $q\Bo$
      and $y^*\Bo y$ has trace zero.
    Regrouping gives a traveling product in
      $\Omega_1=\{y,y^*\}\cup qw\NvNo_1w^*\cup yw\NvNo_1w^*\cup
      w\NvNo_1w^*y^*\cup yw\NvNo_1w^*y^*$ and $\Bo$.
    Now $\Omega_1\subset\NvN_1$,
      $\lspan\Lambda(\{a\}\cup\Ao\cup a\Ao,\{b\})$ is a dense
      $*$--subalgebra of $\NvN_1$
      and $\tau(z)=\tau(qz)=0$ $\forall\;z\in\Omega_1$, so that as above,
      each $z\in\Omega_1$ is the s.o.--limit of a bounded sequence in
      $\lspan(\Lambda(\{a\}\cup\Ao\cup a\Ao,\{b\})\backslash\{1,b\})$.
    So it suffices to show that a nontrivial alternating product in
      $\lspan(\Lambda(\{a\}\cup\Ao\cup a\Ao,\{b\})\backslash\{1,b\})$ and $\Bo$
      has trace zero.
    Regrouping gives a nontrivial alternating product in
      $\{a\}\cup\Ao\cup a\Ao$ and $\{b\}\cup\Bo\cup b\Bo$,
      which by freeness has trace zero.

    Now we prove~(ii).
    Let $\MvN$ be the von Neumann algebra on the left hand side of~(ii),
      and let $\tau$ be its trace.
    We will identify $A$ with $A\otimes1$ and $B$ with $B\otimes1$ as in
      the proof of~(i).
    Let $p$ be a projection in $1\otimes M_2$ (commuting with $A$)
      of trace $\frac12$ and $q$ a projection in $1\otimes L(\Intg_2)$
      (commuting with $B$) of trace $\frac12$.
    Let $\NvN_0=\{p,q\}''$ and let $x,y,w,a,b\in\NvN_0$ be as in the
      proof of~(i).
    Let $u\in1\otimes M_2$ be a partial isometry from $p$ to $1-p$.
    Then
      $$ p\MvN p=\bigl(\{pqp,x^*u\}\cup pA\cup w^*qBw\cup w^*y^*Byw\bigr)'', $$
      and we shall show that $x^*u$ is a Haar unitary ({\it i.e\.} a unitary
      such that $(x^*u)^n$ has trace zero $\forall\;n\in\Intg\backslash\{0\}$)
      and that
      $\bigl\{\{pqp\},\{x^*u\},pA,w^*qBw,w^*y^*Byw\bigr\}$ is $*$--free in
      $p\MvN p$.
    This will in turn prove~(ii).
    For $n>0$, $r=(x^*u)^n$ is a nontrivial alternating product in $\{x^*\}$
      and $\{u\}$, and $x^*$ is the s.o.--limit of a bounded sequence in
      $\lspan(\Lambda(\{a\},\{b\})\backslash\{1,a\})$, so to show $\tau(r)=0$
      it suffices to show that a nontrivial alternating product in
      $\lspan(\Lambda(\{a\},\{b\})\backslash\{1,a\})$ and $\{u\}$ has trace
zero.
    Regrouping gives a nontrivial alternating product in $\{a,u\}$
      and $\{b\}$, which by freeness has trace zero.
    Hence we have shown that $x^*u$ is a Haar unitary in $p\MvN p$.

    Now we show that $x^*u$ and $pqp$ are $*$--free in $p\MvN p$.
    Let $g_k$ $(k\ge1)$ be as in the proof of~(i).
    It suffices to show that a nontrivial alternating product in
      $\{(x^*u)^n\mid n\in\Intg\backslash\{0\}\}$ and
      $\{g_k\mid k\ge1\}$ has trace zero.
    Regrouping gives an alternating product in $\Omega_0$ and $\{u,u^*\}$,
      where $\Omega_0$ is as in the proof of~(i), which, proceeding as we
      did above, we see has trace zero.
    Similarly, we can show that letting $\NvNw_0=\{pqp,x^*u\}''$,
      $\{\NvNw,pA\}$ is free in $p\MvN p$, and that
      letting $\NvNw_1=(\NvNw_0\cup A)''$,
      $\{w^*\NvNw_1w,qB,y^*By\}$ is free in $q\MvN q$, thus proving~(ii).

    To prove~(iii), let $p$ and $u$ in $1\otimes M_2$ commuting
      with $A$ be as above, let $q\in1\otimes M_2$ commuting with $B$ be
      a projection of trace $\frac12$ and $v\in1\otimes M_2$ commuting
      with $B$ a partial isometry from $q$ to $1-q$.
    Let $x,y,w\in\NvN_0=\{p,q\}''$ be as above.
    Then we similarly show that $x^*u$ and $y^*v$ are Haar unitaries
      and that
      $\bigl\{\{pqp\},\{x^*u\},pA,\{w^*y^*vw\},w^*qBw,\bigr\}$ is $*$--free in
      $p\MvN p$, (and notice that these taken together generate $p\MvN p$),
      which proves~(iii).
  \QED

  \proclaim{Corollary 3.6}
    Let $R$ and $\Rw$ be copies of the hyperfinite II$_1$--factor.
    Then
      $$ R*\Rw\cong R*L(\Intg), $$
      with an isomorphism which when restricted is the identity map
      from $R$ to $R$.
  \endproclaim
  \demo{Proof}
    Write $R=(pRp)\otimes M_2$ and $\Rw=(\pw\Rw\pw)\otimes M_2$, where
      $p$ and $\pw$ are projections of trace $\frac12$ in $R$ and
      respectively $\Rw$.
    Then by~(iii) and the proof of~(iii),
      $$ p(R*\Rw)p\cong(pRp)*(\pw\Rw\pw)*L(\freeF_3), $$
      and the isomorphism when restricted to $pRp\subset p(R*\Rw)p$ is
      the identity map from $pRp$ to $pRp$.
    Similarly, writing also $L(\Intg)\cong L(\Intg)\otimes L(\Intg_2)$,
      we have from~(ii) and the proof of~(ii) that
      $$ p(R*L(\Intg))p\cong(pRp)*L(\freeF_4), $$
      and the isomorphism, when restricted to $pRp\subset p(R*L(\Intg))p$, is
      the identity map from $pRp$ to $pRp$.
    Considering the isomorphism~(3), we get an isomorphism from
      $p(R*\Rw)p$ to $p(R*L(\Intg))p$ which when restricted is the identity
      map on $pRp$.
    Now tensor with $M_2$.
  \QED

\noindent{\bf \S4\. The addition formula for free products.}

  \proclaim{Theorem 4.1}
    $L(\freeF_r)*L(\freeF_{r'})=L(\freeF_{r+r'})$ for $1<r,r'\le\infty$.
  \endproclaim
  \demo{Proof}
    (Please see the comments at the end of the introduction.)
    In a W$^*$--probability space $(\MvN,\tau)$ where $\tau$ is a trace,
      let $R$ and $\Rw$ be copies of the hyperfinite II$_1$--factor
      and let $\nu=\{X^t\mid t\in T\}$ be a semicircular family such that
      $\{R,\Rw,\nu\}$ is free.
    Let
      $$ \align
        L(\freeF_r)=\Ac & =\bigl(R\cup\{p_sX^sp_s\mid s\in S\}\bigl)'', \\
        L(\freeF_{r'})=\Bc & =\bigl(\Rw\cup\{q_sX^sq_s\mid s\in S'\}\bigl)'',
      \endalign $$
      where $S$ and $S'$ are disjoint subsets of $T$, $p_s\in R$,
      $q_s\in\Rw$ are projections and where $1+\sum_{s\in S}\tau(p_s)^2=r$,
      $1+\sum_{s\in S'}\tau(q_s)^2=r'$.
    Then $\Ac$ and $\Bc$ are free in $(\MvN,\tau)$, so
      $$ L(\freeF_r)*L(\freeF_{r'})\cong
        \NvN=\bigl(R\cup\Rw\cup\{p_sX^sp_s\mid s\in S\}
        \cup\{q_sX^sq_s\mid s\in S'\}\bigl)''. $$
    By Corollary~3.6, there exists a semicircular element
      $Y\in\NvN_0=(R\cup\Rw)''$ such that $R$ and $\{Y\}$ are free and
      together they generate $\NvN_0$.
    Moreover, for $s\in S'$ let $U_s\in\NvN_0$ be a unitary such that
      $U_sq_sU_s^*=f_s\in R$.
    Then
      $$ \NvN=\bigl(R\cup\{Y\}\cup\{p_sX^sp_s\mid s\in S\}\cup
        \{f_s(U_sX^sU_s^*)f_s\mid s\in S'\}\bigl)''. $$
    To prove the theorem, it suffices to observe that
      $\{R,\{Y\},(\{X^s\})_{s\in S},(\{U_sX^sU_s^*\})_{s\in S'}\}$
      is free in $\MvN$.
  \QED

  Let us recall~\cite{4}
    that the fundamental group of a II$_1$--factor $\MvN$ is defined
    to be the set of positive real numbers $\gamma$ such that
    $\MvN_{\gamma}\cong\MvN$.
  Murray and von Neumann~\cite{4}
    showed that the fundamental group of the hyperfinite
    II$_1$--factor is $\Reals_+$, and recently R\u{a}dulescu~\cite{5} has shown
    that the fundamental group of $L(\freeF_\infty)$ is also $\Reals_+$.
  A\. Connes~\cite{1}
    has shown that the fundamental group of $L(G)$ where $G$
    is a group with property~T of Kazhdan must be countable, but
    no other examples are known for fundamental groups of II$_1$--factors.

  Equation~(2) shows that the isomorphism question for (interpolated)
    free group factors is equivalent to the fundamental group question.
  Combined with the addition formula for free products, we now see that we
    must have one of two extremes.
  \proclaim{Corollary 4.2}
    We must have either
      \roster
      \item"(I)" $L(\freeF_r)\cong L(\freeF_{r'})$ for all $1<r,r'<\infty$ and
        the fundamental group of $L(\freeF_r)$ is $\Reals_+$ for all
        $1<r<\infty$,
      \item"or (II)" $L(\freeF_r)\not\cong L(\freeF_{r'})$
        for all $1<r<r'<\infty$ and the fundamental group of
        $L(\freeF_r)$ is $\{1\}$ for all $1<r<\infty$.
      \endroster
  \endproclaim
  \demo{Proof}
    Using formulas~(1) and~(2) we can show that if
      $L(\freeF_r)=L(\freeF_{r'})$ for some $r\neq r'$, then we have
      $L(\freeF_r)=L(\freeF_{r''})$ for $r''$
      in some open interval, hence that the fundamental
      group of $L(\freeF_r)$ contains an open interval, thus is all of
      $\Reals_+$.
  \QED

\noindent{\bf Acknowledgements.}

  I would like to thank Dan Voiculescu, my advisor, for helpful discussions
    and for suggesting I look at free products such as $M_2(\Cpx)*M_2(\Cpx)$.

\end{document}